\documentclass[twoside]{article}
\usepackage[utf8]{inputenc}
\usepackage{amsmath}
\usepackage{amsfonts}
\usepackage{amssymb}
\usepackage{graphicx}
\usepackage{color}
\usepackage{float}
\usepackage[left=2cm,right=2cm,top=2cm,bottom=2cm]{geometry}
\usepackage{fdsymbol}
\newenvironment{Proof}{\noindent{\sc Proof.}}{\qed}
\newtheorem{theorem}{Theorem}
\newtheorem{lemma}{Lemma}[section]

\newcommand{\qed}{\hfill$\Box$\par\medskip}

\def\bhag#1{\noindent
\setcounter{equation}{0}
\section{#1}
}

\def\CC{{\mathbb C}}

\def\PPI{{{\rm I}\kern-1pt\Pi}}

\def\b #1;{{\bf #1}}

\def\be{\begin{equation}}
\def\ee{\end{equation}}
\def\bea{\begin{eqnarray}}
\def\eea{\end{eqnarray}}

\def\donchitre#1#2{\vskip 6.5cm\noindent
\parbox[t]{1in}{\special{eps:#1.eps x=6.5cm y=5.5cm}}
\hbox to 7cm{}\parbox[t]{0.0cm}{\special{eps:#2.eps x=6.5cm y=5.5cm}}}

\def\CC{{\mathbb C}}

\title{On Marcinkiewicz-Zygmund inequalities and $A_p$-weights for $L$-shape arcs}
\author{Charles K. Chui$^{}$\thanks{Department of Statistics, Stanford University, CA 94305;
 \textsf{email:} ckchui@stanford.edu. The research of this author was partially supported by the Hong Kong Research Council, under Grant $\sharp$ 12303218}
\ and
Lefan Zhong$^{}$\thanks{\textsf{email:} lefan$\_{\rm z}$@yahoo.com}}
\date{}
\begin{document}

\maketitle
 \begin{abstract}
Let $\Gamma$ be an $L$-shape arc consisting of 2 line segments that meet at an angle different from $\pi$ in the complex $z$-plane $\CC$. Application of the exterior conformal map $\psi$ from $|w| > 1$ onto $\CC\backslash\Gamma$, with $\psi(\infty)= \infty$,  introduces the level curves $\Gamma_n=\{z= \psi(w):|w|=1+{1\over{n+1}}\}$. Let $\psi^*$ denote the continuous extension of $\psi$ from $|w|> 1$ to $|w|\ge 1$, so that any family $\{z_{n,k}: k = 0, 1, \dots, n\}$ of points on $\Gamma$ can be written as $\{z_{n,k} = \psi^*(w_{n,k})\}$, where  $|w_{n,k}|= 1$. 
Let $\omega_n (z)= \Pi^n_{k=0} (z-z_{n,k})$. The main objective of this paper is to show that for $L$-shape arcs, validation of the Marcinkiewicz-Zygmund inequalities is equivalent to that of the totality of the $A_p$-weight conditions of $|\omega_n (z) |$ on $\Gamma_n$ and a mild separation condition of $\{z_{n,k}\}$. Since the Marcinkiewicz-Zygmund inequalities are essential to the study of Lagrange polynomial interpolation of continuous functions at the nodes  $\{z_{n,k}\}$, another objective of this paper is to investigate the behavior of the polynomial interpolants at the Fej\'er points, defined by $\{z_{n,k} = \psi^*(e^{i(2k\pi + \theta)/(n+1)})\}$ for any choice of $\theta$. In this regard, we recall that for the interval [-1, 1], the Fej\'er points $\{z_{n,k} = \psi^*(e^{i(2k+1)\pi/(n+1)})\}$ agree with the Chebyshev points and that the Chebyshev points are most commonly used as nodes for Lagrange polynomial interpolation.  On the other hand, numerical experimentation demonstrates that for a typical open $L$-shape arc $\Gamma$, the Lebesgue constants tend to $\infty$ at the rate of $O((log(n))^2)$, as the polynomial degree $n$ increases, while the $A_{p}$-weight conditions for the Fej\'er points $\{z_{n,k}\}$ do not carry over from [-1, 1] to a truly $L$-shape arc. Further numerical experiments also demonstrate that the least upper bounds of the Marcinkiewicz-Zygmund inequalities for the canonical Lagrange interpolation polynomials at $\{z_{n,k}\}$ seem to grow at the rate of $n^{\beta}$, for some $\beta >0$ that depends on $p >1$.
\end{abstract}
\vskip 10pt
 \begin{center} \bf{Dedicated to Professor Guido Weiss on the occasion of his 90-th Birthday}
 \end{center}

 \bhag{Introduction}
Let $\Pi_n$ denote the space of polynomials with degree at most $n$, and recall the well-known Marcinkiewicz-Zygmund inequalities $\cite{MZ}$:
\begin{equation}
{{c_1}\over{n+1}} \sum^n_{k=0} \Big|P_n \Big(e^{i{{2\pi k}\over{{n+1}}}}\Big)\Big|^p \le {1\over{2\pi}}\int_0^{2\pi} \Big|P_n (e^{i\theta}) \Big|^p  d\theta \le {{c_2}\over {n+1}} \sum^n_{k=0} \Big|P_n \Big(e^{i{{2\pi k}\over{{n+1}}}}\Big)\Big|^p,
\end{equation}
for $1<p< + \infty,   P_n\in \Pi_n$, where $0< c_1 \le c_2$, are constants independent of $n$ and maybe different for other inequalities throughout the discussions in this paper. It is natural to ask if the roots of unity in (1.1) can be replaced by a general family of $\{z_{n,k}: k=0,1,\dots,n\}$ on the unit circle. In this regard, it is shown in our earlier work $\cite{CZ}$, that the validity of the general Marcinkiewicz-Zygmund inequalities on the unit circle:
\begin{equation}
{{c_1}\over {n+1}}\sum_{k=0}^n \big|P_n (z_{n,k}) \big|^p \le {1\over 2\pi}\int_{0}^{2\pi}\big|P_n(e^{i\theta})\big|^p  d\theta \le {{c_2}\over {n+1}} \sum_{k=0 }^n \big|P_n (z_{n,k} ) \big|^p,
\end{equation}
for $1<p< + \infty$, is equivalent to the validity of the $A_p$-weight condition of $\omega_n(z) = \Pi_{k = 0}^n(z - z_{n,k})$:
 \begin{equation}
\sup_{t_0<t_1 } \Big \{ {1\over{t_1-t_0}} \int_{t_0}^{t_1} | \omega_n ((1+{{1}\over {n+1}}) e^{it} ) |^p dt \Big \}^{1\over p} \Big \{{1\over {t_1-t_0}} \int_{t_0}^{t_1} |\omega_n ((1+{{1}\over {n+1}}) e^{it}) |^{-q} dt \Big\}^{1\over q}\le c_2,
\end{equation}
where ${1\over{p}} + {1\over{q}} = 1$, together with the separation condition:

\begin{equation}
 \min_{0\le k<n} |z_{n,k+1} - z_{n,k} | \ge {{c_1}\over{n+1}}.
\end{equation}
This establishes a close relationship of the Marcinkiewicz-Zygmund inequalities with the $A_p$-weights condition, introduced by Muckenhoupt $\cite{M}$. To the best our knowledge, the intimate relationship between the two seemingly different research topics: Marcinkiewicz-Zygmund ($M$-$Z$) inequalities and Muckenhoupt $A_p$-weights ($M$-$A_p$) condition, was initiated in our earlier work $\cite{CZ}$; and unfortunately, we do not know of any serious follow-up research development in this direction in the published literature. The importance for a family of points $z_{n,k}$, $k=0,\cdots,n$, to satisfy the $M$-$A_p$ condition is that it is necessary condition for this family to replace the $(n+1)$-$th$ roots unity in the $M$-$Z$ inequalities (1.1).
\vskip 5pt
From the Approximation Theory point of view, the Marcinkiewicz-Zygmund inequalities (1.2) for a general family of points on the unit circle are also equivalent to the optimal order $O({1\over n})$ of $H^p$ approximation by the Lagrange interpolation polynomials $L_n (f,z)$ at $\{z_{n,k}\}$  of $f(z)\in C(|z|\le 1)$, with $f^\prime(z)\in H^p(|z|<1)$, as shown in our earlier work $\cite{CZ}$. In other developments, different weights of the Marcinkiewicz-Zygmund inequalities for $z_{n,k}\in [-1, 1]$ and $z_{n,k}\in(-\infty, \infty)$ are studied in $\cite{L}$ by using different nodes $\{z_{n,k}\}$ including Jacobi zeros; with Chebyshev nodes on [-1, 1] extended to [-1, 1]x[-1, 1] in the paper $\cite{Xu}$; and the inequalities in (1.2) generalized to $p = 1$ and $\infty$ in $\cite{OCS}$, by using $m = (1+\epsilon)n$ nodes $\{z_{m,k}: k=0\ldots,m\}$ for some $\epsilon > 0$ with $|z_{m,k}| = 1$. In addition, the results in $\cite{OCS}$ are also extended from $|z| = 1$ to the unit sphere of dimension higher than 1 in  $\cite{W}$. However, all of the considerations mentioned above are for such regular point-sets as circles/spheres, lines/intervals, and unit square. 
\vskip 5pt

 Let $U$ denote the unit disk $|w| < 1$ with closure $\bar U$ and consider a one-dimensional compact and connected rectifiable curve $\gamma$ in the complex plane $\CC$ that does not cross itself, such as a simple closed curve or a simple open arc. Then the conformal map $\psi$ from $\CC \backslash \bar U$ to  complement of $\gamma$ with $\psi(\infty)= \infty$ can be extended to a continuous function  $\psi^*$ on $|w| \ge 1$. Let $\gamma_n$ be the level curve $\{\psi(w): |w|= 1+{1\over {n+1}}\}$, and denote by $dist(z, \gamma_n)$ the (Euclidean) distance from $z$ to set $\gamma_n$. A family $\{z_{n,k}:k=0,1,\dots,n\}$ of nodes that lie on the set  $\gamma$  will be called a Marcinkiewicz-Zygmund ($M$-$Z$) family on $\gamma$, if any polynomial $P_n \in \Pi_n$ satisfies: 
\begin{equation} \label{ineq:MZ}
c_1 \sum^n_{k=0} |P_n (z_{n,k} ) |^p dist(z_{n,k},\gamma_n )    \le \int_\gamma |P_n (z) |^p |dz| \le  c_2 \sum^n_{k=0} |P_n (z_{n,k} ) |^p dist(z_{n,k},\gamma_n ).
\end{equation}
This is a generalization of (1.2) from unit circle to $\gamma$.

To generalize to $M$-$A_p$-weight condition (1.3), we introduce the notion of Muckenhoupt $A_p$-weight constants: 
\begin{equation} \label{ineq:Mn}
M_n = { sup}_{I \subset \gamma_n }    { \Bigg({1\over {|I|}} \int_I |\omega_n (\zeta)|^p |d\zeta| \Bigg)^{{1\over p}} \ \Bigg({1 \over {|I|}} \int_I |\omega_n(\zeta)| ^ {-q} |d\zeta| \Bigg)^{1\over {q}} },
\end{equation}
corresponding to the family $\{z_{n,k}:k=0,1,\dots,n\} \subset \gamma$ of points. The family $\{z_{n,k}\}$ is said to satisfy the $M$-$A_p$ condition, if the sequence  $\{M_n\}$ of Muckenhoupt $A_p$-weight constants are bounded by some constant $M>0$.

Finally, the separation condition (1.4) is generalized to the $c_0$-separation condition, defined by 
\begin{equation} \label{ineq:sep_c0}
\min_{0\le k < j\le n}   {{|z_{n,j}-z_{n,k}|} \over {\min \big(dist(z_{n,j},\gamma_n ),dist(z_{n,k},\gamma_n ) \big)}} \ge c_0>0,
\end{equation}
where $\{z_{n,k}:k=0,1,\dots,n\} \subset \gamma$.

To extend the Marcinkiewicz-Zygmund inequalities (1.1) from the unit circle to a general $\gamma$, it is important to be able to identify $M$-$Z$ families ${\{z_{n,k}}\}$. The common approach is to consider the set  ${\{\psi((1+ 1/m)e^{i(2k\pi + \theta)/(n+1)})\}}$ of Fej\'er points on $\gamma_m$ and let $m$ tend to infinity to arrive at 
\begin{equation} \label{eq:Znk}
\{z_{n,k}=\psi^*(e^{i(2k\pi + \theta)/(n+1)})\}\subset \gamma.
\end{equation}
Indeed, for a simple closed (Jordan) curve in $C^{2+\epsilon}$, it was shown in \cite{AK} that the Fej'er points in (\ref{eq:Znk}) is an $M$-$Z$ family, and this result was extended to piecewise $C^{2+\epsilon}$ Jordan curves with corners (see Lemma 6 in \cite{Z}). However, for open arcs $\gamma$, the situation is more complicated, since by taking the limit of   ${\{\psi((1+ 1/m)e^{i(2k\pi + \theta)/(n+1)})\}}$, it is quite possible that two points of the family $\{z_{n,k}:k=0,1,\ldots,n\}$ in (\ref{eq:Znk}) may become one single point. In fact, $\psi^*$ is not an univalent function on the unit circle. In addition, the separation condition (1.8) is no longer assured for ${z_{n,k}}$ in (\ref{eq:Znk}). In this regard, we mention that the method developed in \cite{ZZ} can be applied to adjust points in (\ref{eq:Znk}) to obtain an $M$-$Z$ family. When $\gamma$ is the interval [-1, 1], the Chebyshev nodes $\{cos({{(2k+1)\pi}\over{2N}}: k=0,1,\ldots,N-1\}$ can be derived from (\ref{eq:Znk}) by $n = 2N - 1$ and $\theta = \pi$ with 2 Fej\'er points mapping to a Chebyshev node.
\vskip 5pt

The motivation of our research in this paper is to establish an analogous equivalence relationship between $M$-$Z$ inequalities and $M$-$A_p$ condition for an open arc. While the task is quite daunting when the arc consists of corners, we have decided to focus our research on establishing the equivalence relationship between $M$-$Z$ inequalities and $M$-$A_p$ condition for an open $L$-shape arc $\Gamma$ which consists of two line segments that meet at an angle different from $\pi$. This happens to be a lucky decision, since we found out from numerical experiments the significant distinction between smooth arcs $\gamma$ without corners and an $L$-shape arc $\Gamma$ in the choice of the family of points $\{z_{n,k}: k=0,\cdots,n\}$ on $\gamma$ and $\Gamma$. We will first show that any family $\{z_{n,k}: k = 0,\cdots, n\} \subset \Gamma$ is an $M$-$Z$ family on  $\Gamma$, if and only if it satisfies both the $M$-$A_p$ condition $M_n\le M$ and the $c_0$-separation condition (\ref{ineq:sep_c0}). 

Another objective of this paper is to study if the family of Fej\'er points $\{z_{n,k}: k = 0,\cdots, n\} \subset \Gamma$ an $M$-$Z$ family on $\Gamma$. In this regard, our numerical results for a typical $L$-shape arc indicate that $M$-$A_p$ condition seems to fail for the family of Fej\'er points, implying that the $M$-$Z$ inequalities (\ref{ineq:MZ}) might not be satisfied by certain polynomials. Such phenomena is very different from smooth open arc $\cite{ZZ}$ (including the interval [-1,1]). It is also different for simple close curves with corners (such as the boundary of a polygon) $\cite{Z}$.

Since an $M$-$Z$ family on $\Gamma$ serves as a useful measuring stick for selecting interpolation nodes for Lagrange polynomial interpolation, the third objective of this paper is to study if the uniform approximation order of  Lagrange polynomial interpolation is controlled by the growth of the Lebesgue constants, defined by
\begin{equation} \label{eq:Ln}
{L_n}  :=  \max_{z\in \Gamma} \sum_{k=0}^n {{|\omega_n (z)|}\over {|\omega_n^\prime (z_{n,k}) } (z-z_{n,k} ) |} .
\end{equation}
Recall that for Lagrange polynomial interpolation at the Chebyshev nodes on [-1, 1], the Lebesgue constants have growth order ${L_n} = O(log(n))$. In $\cite{ZZ}$, this growth order is also established for Lagrange polynomial interpolation at adjusted Fej\'er points on a $C^{2+\epsilon}$ smooth open arc. In the present paper, we will demonstrate with numerical experiments that after the Fej\'er points on a typical $L$-shape arc are adjusted to satisfy separation condition (1.7), the Lebesgue constants ${ L_n}$ for Lagrange polynomial interpolation at $\{z_{n,k}^*\}$ seem to have growth order of $O((log(n))^2)$.

\bhag{Marcinkiewicz-Zygmund inequalities and $A_p$-weights for $L$-shape arcs}
Let $\Gamma$ denote an $L$-shape arc in the complex plane $\CC$ that consists of two 2 line segments which meet at an angle different from $\pi$. That is, the one-dimensional compact and connected set $\gamma$ and its corresponding level curves $\gamma_n$ introduced in the previous section are replaced by $\Gamma$ and $\Gamma_n$, respectively, in the following discussions. This section is devoted to the proof of the following theorem.

\begin{theorem}  
Let $\Gamma$ be an $L$-shape arc in $\CC$ and $1<p< +\infty$. Then $\{{z_{n,k}: k=0,1,\dots,n}\} \subset \Gamma$ is an $M$-$Z$ family if and only if it satisfies both of the $M$-$A_p$ and the $c_0$-separation conditions. That is, the following two statements are equivalent.

\begin{itemize}
  \item[(i)]
	All polynomials $P_n \in \Pi_n$ satisfy:
	
\begin{equation} \label{ineq:MZ2}
c_1 \sum^n_{k=0} |P_n (z_{n,k} ) |^p dist(z_{n,k},\Gamma_n )    \le \int_\Gamma |P_n (z) |^p |dz| \le  c_2 \sum^n_{k=0} |P_n (z_{n,k} ) |^p dist(z_{n,k},\Gamma_n ).
\end{equation}
	
 \item[(ii)] 
	There exist constants $M>0$ and $c_0 > 0$, such that 
\begin{equation} \label{ineq:Mn2}
M_n := {sup}_{I \subset \Gamma_n }    { \Bigg({1\over {|I|}} \int_I |\omega_n (\zeta)|^p |d\zeta| \Bigg)^{{1\over p}} \ \Bigg({1 \over {|I|}} \int_I |\omega_n(\zeta)| ^ {-q} |d\zeta| \Bigg)^{1\over {q}} }\le M,
\end{equation}	
and 
\begin{equation} \label{ineq:sep_c02}
\min_{0\le k < j\le n}   {{|z_{n,j}-z_{n,k}|} \over {\min \big(dist(z_{n,j},\Gamma_n ),dist(z_{n,k},\Gamma_n ) \big)}} \ge c_0>0.
\end{equation}
\end{itemize}
\end{theorem}

We first remark that if we re-arrange the points $\{z_{n,k}:k=0,1,\dots,n\}$ along $\Gamma$ in an increasing order, then it is easy to see that the condition (\ref{ineq:sep_c02}) is equivalent to 
\begin{equation} \label{ineq:sep_c1}
\min_{0\le k<n}  |z_{n,k+1}-z_{n,k} |  \ge c_1 \min \big( dist(z_{n,k+1},\Gamma_n ),dist(z_{n,k},\Gamma_n ) \big).
\end{equation}
Here, for convenience, we have used the same notation for the re-arranged points. To establish the theorem, we need the following three lemmas.
\begin{lemma}
For any polynomial $P_n \in \Pi_n$, the integrals of the $L_p$-norm over $\Gamma$ and $\Gamma_n$ are equivalent to the $L_p$ norm, in that 
\begin{equation} \label{ineq:Ben}
c_1 \int_{\Gamma_n} |P_n (z) |^p |dz|\le \int_\Gamma |P_n (z) |^p |dz| \le c_2 \int_{\Gamma n} |P_n (z) |^p |dz|. 
\end{equation} 
\end{lemma}                                                      
The above lemma is a re-formulation of Lemma 10 in $\cite{D}$.

Next let $D_n$ be the domain with boundary $\Gamma_n$ and denote by $E^p (D_n)$ the Hardy space on $D_n$ , i.e. the space of functions $f$ analytic in $D_n$ with $f|_{\Gamma_n}\in L^p (\Gamma_n)$.

\begin{lemma}
Under the $c_0$-separation condition in (\ref{ineq:sep_c02}),
\begin{equation} \label{ineq:carlson}
\sum^n_{k=0} |f(z_{n,k} ) |^p  dist(z_{n,k}, \Gamma_n) \le  c_2 \int_{\Gamma_n} |f(z) |^p |dz|, 
\end{equation}
for all $f\in E^p (D_n )$.
\end{lemma}

\begin{Proof} Following the proof of Lemma 5.2 in $\cite{ZZ}$, we consider the disjoint union of circles
$$V_n= \bigcup_{k=0}^n \{z: |z-z_{n,k} | = {{c_0}\over {2}} dist(z_{n,k},\Gamma_n)\}.$$
In view of (\ref{ineq:sep_c02}), we note that $V_n$ is a disjoint union of $n+1$ circles with centers at $\{z_{n,k}: k=0,1,\dots,n\}$ in $\Gamma$. Hence, since $|f(z) |^p$ is subharmonic in $D_n$, we have
$$ \sum_{k=0}^n |f(z_{n,k} ) |^p  dist(z_{n,k},  \Gamma_n) \le {{2}\over {c_0}}   {1\over {2\pi}} \int_{V_n} |f(z) |^p |dz| \le c_2 \int_{\Gamma_n} |f(z) |^p |dz|.$$
The second inequality holds due to the  Carleson measure behavior of $V_n$, in that the arc length of $V_n$ inside any disk $U(z, R)$, with center at $z$ and radius $R$ can be controlled by $c_2 R$.
\end{Proof}

\begin{lemma}
Let $\{z_{n,k}\}\subset\Gamma$ satisfy  (2.3) and $\{a_{n,k}:k=0,1,\ldots,n\}\subset\CC$ be arbitrarily given. Then there exists an $f \in E^p (D_n )$, such that 
$$f(z_{n,k}  )=a_{n,k},$$
for   $k=0,1,\ldots,n$, and 
\begin{equation}
\int_{\Gamma_n} |f(z) |^p |dz| \le c_{2}  \sum_{k=0}^n |a_{n,k} |^p  dist(z_{n,k}, \Gamma_n ).
\end{equation}
\end{lemma}
\begin{Proof}
Let $\delta_z$ denote the linear functional $\delta_z f = f(z)$ for $f\in E^p(D_n)$. As in the proof of Lemma 3 in $\cite{Z}$, since the distribution
$$\sum_{k=0}^n \delta_{z_{n,k}} dist(z_{n,k},\Gamma_n)$$
is a Carleson measure in $D_n$ and $\{z_{n,k}\}$ satisfies the separation condition (2.3), the proof of Lemma 2.3 is complete by applying Theorem 1.1 in chapter VII of $\cite{G}$. 
\end{Proof}

We are now ready to prove Theorem 1. Our proof follows the idea of the derivation in our earlier work $\cite{CZ}$.
\vskip 5pt

To see that $(ii)$ implies $(i)$, we simply observe that under the condition in (2.3), it follows from Lemma 2.1 that
$$\sum_{k=0}^n |P_n (z_{n,k} ) |^p dist(z_{n,k},\Gamma_n ) \le c_2 \int_{\Gamma_n} |P_n (z) |^p |dz| \le c_2  \int_{\Gamma} |P_n (z) |^p |dz|,$$
which establishes the lower bound in (2.1). To prove the upper bound in (2.1), we consider the function $f$ in Lemma 2.3 that interpolates the data $\{P_n  (z_n,k):k=0,1,\dots,n\}$, i.e. $f(z_{n,k} ) = P_n (z_{n,k}),  k=0,1,\ldots,n$, so that Lemma 2.3 yields  
\begin{equation}
\int_{\Gamma_n} |f(z) |^p |dz| \le c_{2}  \sum_{k=0}^n |P_n (z_{n,k})|^p  dist(z_{n,k}, \Gamma_n ).
\end{equation}
Next, we use the Cauchy singular integral operator, defined by
\begin{equation}
(Sg)(z) := {1\over{2 \pi i}}\int_{\Gamma_n} g(\zeta)  {d\zeta \over {\zeta - z}}, 
\end{equation}
to formulate 
\begin{equation}
f(z)-P_n (z)={{\omega_n (z)}\over {2\pi i}}  \int_{\Gamma_n}{{f(\zeta)}\over{\omega_n (\zeta)}} {{ d\zeta} \over {\zeta - z}} = \omega_n(z) (S\big ({f\over {\omega_n}}\big))(z).
\end{equation}                    
Therefore, it follows from the Minkowski Inequality that

\begin{equation}	
\big\{\int_{\Gamma_n} |P_n (z) |^p |dz|\big\}^{1\over p} \le \big\{\int_{\Gamma_n} |f(z) |^p |dz|\big\}^{1\over p}  + \big\{\int_{\Gamma_n} |(S\big({f\over {\omega_n}}\big)) (z) \big|^p |\omega_n  (z) |^p |dz|\big\}^{1\over p}. 
\end{equation}
Now, since $|\omega_n (z)|$ satisfies the $A_p$-condition in (2.2), the Cauchy singular integral operator $S$ is bounded in $L^p (\Gamma_{n},|\omega_n |^p )$, so that
$$\int_{\Gamma_n} \big|S \Big({f\over {\omega_n}}\Big)(z) \big|^p | \omega_n (z) |^p |dz| \le c_{2} \int_{\Gamma_n} \big| \Big({f\over {\omega_n}}\Big) (z) \big|^p |\omega_n  (z) |^p |dz| = c_{2} \int_{\Gamma_n} |f(z) |^p |dz|. $$ 
Hence, it follows from Lemma 2.1, together with (2.8) and (2.11), that the upper bound in (2.1) also holds.
\vskip 5pt

To prove that $(i)$ implies $(ii)$, 
we first show that the Marcinkiewicz-Zygmund inequalities imply the separation condition (2.4).
Let $\Psi_n$  be the conformal map from $U$ onto $D_n$, with $\Psi_n(0) =  z_{n,0}$ and denote $w_{n,k} = \Psi_n ^{-1}(z_{n,k})$. In view of (1.3), we have, by applying the Koebe Distortion Theorem, that
$${1\over 4} dist(z_{n,k},\Gamma_n )\le |\Psi_n^\prime (w_{n,k} )|\Big(1-|w_{n,k} |^2 \Big) \le 4 dist(z_{n,k},\Gamma_n)$$
Let $P_{n,0} \in \Pi_n$ be defined by $P_{n,0} (z_{n,0} ) = 1$ and $P_{n,0} (z_{n,k} ) = 0$ for $k>0$, and consider the Blaschke product 
$$B_n (w)= \prod_{k>0} {{w-w_{n,k}}\over {1-\overline{w_{n,k}}   w} } .$$
Since $\Psi_n^{\prime} (w) \neq 0$ in $U$, the function ${{P_{n,0} (\Psi_n (w))} \over {B_n (w) }} \Psi_n^\prime (w)^{1/p}$   is analytic in $U$, so that

\begin{eqnarray*}
 \Big|{{P_{n,0} (\Psi_n (w_{n,0}))} \over { {B_n (w_{n,0} ) }}} \Psi_n^\prime (w_{n,0} )^{1\over p} \Big| &=& \Big|{1\over {2\pi}} \int_{|w|=1} {{P_{n,0} (\Psi_n (w))} \over {B_n (w) }} \Psi_n^\prime (w)^{1\over p}  {{dw} \over {(w-w_{n,0})}}\big|  \\
&\le& \Big({1\over 2\pi} \int_{|w|=1} \Big|{{P_{n,0} (\Psi_n (w))} \over {B_n (w)}} \Psi_n^\prime (w)^{1\over p} \Big|^p |dw|  \Big)^{1\over p} \Big({1\over {2\pi}} \int_{|w|=1} \Big|{1\over {w-w_{n,0} }}\Big|^q |dw| \Big)^{1\over q}\\
&=& \Big( {1\over {2\pi}} \int_{\Gamma_n} |P_{n,0} (z)|^p |dz|\Big)^{1\over p} \Big({1\over {2\pi}} \int_{|w|=1} \Big|{1\over {w-w_{n,0} }}\Big|^q |dw| \Big)^{1\over q}\\ 
&\le& c_{2} |P(z_{n,0} ) | (dist(z_{n,0},\Gamma_n ) )^{1\over p} (1-|w_{n,0} |^2 )^{{1-q} \over {q}}
\end{eqnarray*}
i.e.
$${1\over {|B_n (w_{n,0}) |}} \le c_{2} \Big |{ {dist(z_{n,0},\Gamma_n) }\over {\Psi_n^\prime (w_{n,0} )(1-|w_{n,0} |^2 )}} \Big |^{1\over p} \le 4^{1\over p} c_{2}.$$
Thus, we have
\begin{equation}
\Big| \prod_{k>1} {{(w_{n,0}-w_{n,k})}\over {(1-{\overline{w_{n,k}}}  w_{n,0})}}\Big| \ge (4^{1\over p} c_{2})^{-1} > 0.
\end{equation}
In particular, it follows, for $w_{n,0}=0$, that
$$|w_{n,1} | \ge {(4}^{1\over p} c_{2})^{-1}.$$
Therefore, application of the Koebe Distortion Theorem yields
$$\big| \Psi_n (w_{n,1} ) - \Psi_n (w_{n,0} )  \big| \ge {1\over 4} \big |\Psi_n^\prime (w_{n,0}) \big||w_{n,1} | / (1+ |w_{n,1}|)^2;$$ 
that is, 
$${{|z_{n,1}-z_{n,0} |} \over {dist(z_{n,0},\Gamma_n )}} \ge c_{1}$$
which is (2.4) for $k = 0$. The same argument shows that (2.4) holds for any $k \le n$. In fact, (2.12) can be generalized to 
\begin{equation}
min_{j \neq k}\Big| \prod_{k>1} {{(w_{n,j}-w_{n,k})}\over {(1-{\overline{w_{n,k}}}  w_{n,j})}}\Big| \ge c_1 > 0
\end{equation}
In view of (2.13), we have, for $h \in H^p (U)$,
$$\sum_{k=0}^n |h(w_{n,k} ) |^p \big(1 - |w_{n,k}|^2\big) \le c_{2} \int_{|w|=1}  |h(w) |^p |dw|.$$
Hence, for any $H \in E^p (D_n)$, by setting $h(w) = H(\Psi_n (w) ) (\Psi_n^\prime (w) )^{1\over p}$, we may conclude that
\begin{eqnarray*}
\sum_{k=0}^n |H(z_{n,k} ) |^p dist(z_{n,k},\Gamma_n ) &=&  \sum_{k=0}^n |h(w_{n,k} ) |^p {{dist(z_{n,k},\Gamma_n )} \over {|\Psi_n^\prime (w_{n,k} )|}}\\
 &\le& 4 \sum_{k=0}^n |h(w_{n,k} ) |^p (1-|w_{n,k}|^2 ) \le 4 c_{2}  \int_{|w|=1}   |h(w) |^p |dw|;
\end{eqnarray*}
or equivalently,
\begin{equation}
\sum_{k=0}^n |H(z_{n,k} ) |^p dist(z_{n,k},\Gamma_n ) \le  4 c_{2} \int_{\Gamma_n} |H(z) |^p |dz|.
\end{equation}

For any $f \in L^p (\Gamma_n )$, consider $g(z)  = f(z) \omega_n (z)$, and let $P_n \in \Pi_n$, such that $P_n (z_{n,k} ) = (Sg)(z_{n,k})$.  Then for $z \in D_n$, we have
$$ (Sg)(z) = {{\omega_n(z)}\over {2 \pi i}} \int_{\Gamma_n} {{g(\zeta)} \over {\omega_n(\zeta)}} {{d \zeta}\over  {\zeta - z}} + P_n (z).$$
Note that if $g \in E^p (D_n )$ , then this is the same as (2.10), and if $g(z) = {1\over{(z-z_{n,0} )^m}}$, with $m > 0$, then $(Sg)(z) = 0$ so that the formula also holds. In any case, we have
$$(Sg)(z) = \omega_n (z) (Sf)(z) + P_n (z)$$
and
\begin{eqnarray*}
\big(\int_{\Gamma_n} |(Sf)(z)|^p |\omega_n (z) |^p |dz|\big)^{1 \over p} &\le& \big(\int_{\Gamma_n} |(Sg)(z) |^p |dz|\big)^{1 \over p} + \big(\int_{\Gamma_n} |P_n (z) |^p |dz|\big)^{1 \over p} \\  
&\le& \big(\int_{\Gamma_n} |(Sg)(z) |^p |dz|\big)^{1 \over p} + {c_2} \big(\sum_{k=0}^n |P(z_{n,k} ) |^p dist(z_{n,k},\Gamma_n)\big)^{1 \over p}\\
&=& \big(\int_{\Gamma_n} |(Sg)(z) |^p |dz|\big)^{1 \over p} + {c_2} \big(\sum_{k=0}^n |(Sg)(z_{n,k} ) |^p dist(z_{n,k},\Gamma_n)\big)^{1 \over p}.
\end{eqnarray*}

It follows from (2.14), for $H(z) = (Sg)(z)$, that
\begin{eqnarray*}
\big(\int_{\Gamma_n} |(Sf)(z)|^p |\omega_n (z) |^p |dz|\big)^{1 \over p}
&\le& c_2\big(\int_{\Gamma_n} |(Sg)(z) |^p |dz|\big)^{1 \over p}.  
\end{eqnarray*}

Since $S$ is a bounded operator on $L^p(\Gamma_n)$, we have

\begin{eqnarray*}
\big(\int_{\Gamma_n} |(Sf)(z)|^p |\omega_n (z) |^p |dz|\big)^{1 \over p}
&\le& c_2\big(\int_{\Gamma_n} |g(z) |^p |dz|\big)^{1 \over p}\\
&=& c_{2} \big(\int_{\Gamma_n} |f(\zeta) |^p |\omega_n (\zeta)|^p |d\zeta|\big)^{1 \over p}.
\end{eqnarray*}
This shows that the Cauchy Singular Integral operator in uniformly bounded in $L^p (\Gamma_n,|\omega_n  |^p)$. Thus, $|\omega_n|$ are $A_p$-weights that satisfy (2.2) (see Theorem 4.15 in \cite{BK}).
\vskip 25pt

\bhag{Divergent phenomena for $L$-shape arcs from numerical experimentation}
 As a consequence of Theorem 1, it is clear that the product $|\omega_n (z)|$ plays an important role in the study of the Marcinkiewicz-Zygmund inequalities. This section on numerical experimentation is devoted to the investigation of the key properties of $|\omega_n (z)|$, particularly for $L$-shape arcs. In particular, we will calculate the Muckenhoupt $A_p$-weight constants $M_n$ introduced in (2.2) on the level curves $\Gamma_n$  for a typical $L$-shape arc.
We will also calculate the Lebesgue constants and compare their rate of growth with  $log(n)$.

\subsection{A typical $L$-shape arc and its associated exterior conformal map}
Consider the conformal map 
\begin{equation}
\psi_\alpha (w) =\Big(w-{1\over w}\Big) \Big({{w-1}\over {w+1}}\Big)^{\alpha}
\end{equation}
from $|w| > 1$ onto $\CC\backslash\Gamma$, where $\Gamma$
is an $L$-shape arc in the complex plane, consisting of 2 line segments that meet at an angle of $(1-\alpha)\pi$, $0\le \alpha < 1$.
The typical $L$-shape arc $\Gamma$ for our numerical experimentation is the union of two 2 line segments, namely: $[0,{27}^{1\over 4} e^{{{3\pi}\over 4} i}] \cup [0,{27}^{1\over 4} e^{{{5\pi}\over 4} i}]$ with the conformal map
\begin{equation}
\psi (w) =\Big(w-{1\over w}\Big) \Big({{w-1}\over {w+1}}\Big)^{1\over 2}
\end{equation}
The extension of $\psi^*$ of $\psi$ to the unit circle has the properties that $\psi^*(1) = \psi^*(-1) = 0$, and that the end-points of $\Gamma$ are $\psi^*(e^{i {2\over 3} \pi} ) = 27^{1\over 4} e^{i {{3\pi}\over 4}}$ and
$\psi^*(e^{i {4\over 3} \pi} ) = 27^{1\over 4} e^{i {{5\pi}\over 4}}$, as shown in Figure 1, where the level curve $\Gamma_n$, with $n = 8$, is also displayed. 

\begin{figure}[h]
\centering
%\begin{minipage}{0.6\textwidth}
\includegraphics[width=6in,height=5in]{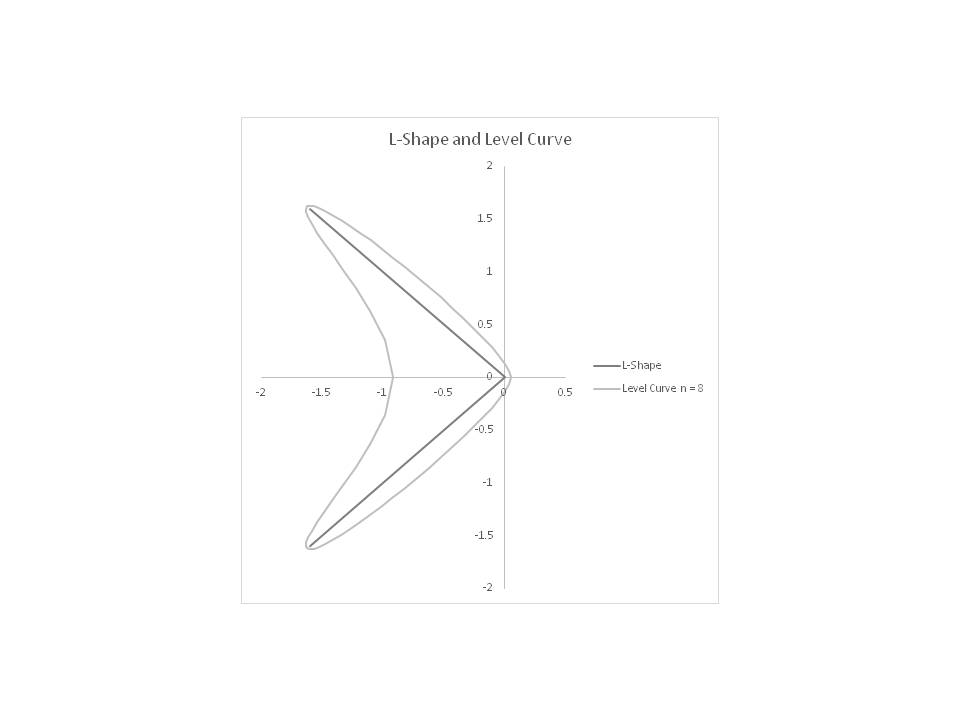} 
%\end{minipage}
\caption{$\Gamma$ and $\Gamma_n$ with $n = 8$}
\end{figure}

\subsection{Minimum values for $|\omega_n(z)|$ on level curves} 

The key to establishing the Marcinkiewicz-Zygmund inequalities in earlier work $\cite{Z}$ and $\cite{ZZ}$ was to prove that
\begin{equation}
{{\max_{z\in \Gamma_n}|\omega_n (z)|}\over {\min_{\zeta\in \Gamma_n}|\omega_n (\zeta)}|} \le c_{2}
\end{equation}                              
on the level curve for $(n+1)-th$ Fej\'er points. Of course (3.3) is much stronger than the $M$-$A_p$-weight condition in (2.2). But as an application of (3.3), along and with (2.3), it is easy prove, using the formula (2.10), that
\begin{equation}
\max_{z\in \Gamma}|P_n (z)| \le c_{2} { log} \ (n) \max_{0\le k \le n}|P_n (z_{n,k} ) |   
\end{equation}                          
for all $P_n\in \Pi_n$. From the same max/min ratio as (3.3), it was shown in [13, Theorem 4] that for smooth open arcs, the order of growth of the Lebesgue constants is $O({ log} \ (n))$ for adjusted Fej\'er points.  This agrees with the order of growth of the Lebesgue constants for the interval $[-1, 1]$ when Chebyshev nodes are considered. Motivated by this result (for smooth open arcs), we decided to investigate the behavior the max/min ratio in (3.3) for $L$-shape open arcs (with corners), by calculating the values of this ratio. To our surprise, numerical experimentation seems to indicate that the ratios increase quite fast for larger values of $n$ (instead of the boundedness in (3.3). In Table 1, we list the values of $\min_{\zeta \in \Gamma_n} |\omega_n (\zeta)|$, $\max_{z \in \Gamma_n} |\omega_n (z)|$ and their ratios, for $n = 2^k, k = 4,5,\ldots, 12$; and in Figure 2, we plot the max/min ratios for even integers $n=10, 12,\ldots, 4096$, and we also noticed that for odd integers $n$, the ratios seem to increase even faster, though they are not addressed in this paper. It also seems that the values of $\min_{\zeta \in \Gamma_n } |\omega_n (\zeta) | \to 0$ as $n \to  \infty$.

\begin{table}[h]
\centering
\begin{tabular}{|r|r|r|r|}
\hline
         \bf $n=2^k$ & \bf $\min_{\zeta\in\Gamma_n} |\omega_n(\zeta)|$ &      \bf  $\max_{z\in\Gamma_n} |\omega_n(z)|$ &  \bf  $\max/\min$ \\
\hline
        16 &    1.09441 &       5.31 &       4.85 \\
\hline
        32 &    0.43913 &       5.59 &      12.73 \\
\hline
        64 &    0.38920 &       5.68 &      14.59 \\
\hline
       128 &    0.43848 &       6.14 &      14.01 \\
\hline
       256 &    0.16583 &       6.41 &      38.64 \\
\hline
       512 &    0.10636 &       6.62 &      62.27 \\
\hline
      1024 &    0.08630 &       6.78 &      78.54 \\
\hline
      2048 &    0.04561 &       6.90 &     151.21 \\
\hline
      4096 &    0.02644 &       6.98 &     264.17 \\
\hline
\end{tabular} 
\caption{$\min$ and $\max$ on $\Gamma_n$}
\end{table}

\begin{figure} [h]
\centering
%\begin{minipage}{0.6\textwidth}
\includegraphics[width=6in,height=4in]{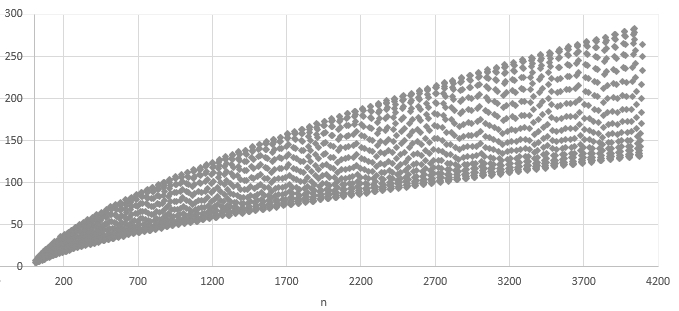} 
%\end{minipage}
\caption{ ${{\max_{z\in \Gamma_n}|\omega_n (z)|}\over {\min_{\zeta\in \Gamma_n}|\omega_n (\zeta)}|}$ on $\Gamma_n$ for even $n = 10, 12,\ldots, 4096$}
\end{figure}

\clearpage
\newpage

\subsection{Lebesgue constants}

Since  (3.3) does not seem to hold for open $L$-shape arcs, perhaps the method introduced in $\cite{ZZ}$ to prove that the growth order $O(log(n))$ of the Lebesgue constants for smooth open arcs cannot be adopted to the study of open $L$-shape arcs. For this reason, we investigate if the Lebesgue constants for open $L$-shape arcs remain to be of the same order as smooth arcs. To carry out this task, we adjust the Fej\'er's points that satisfy the separation condition (2.3) as was done for smooth open arcs in  $\cite{ZZ}$, for computing the Lebesgue constants using the formula

$$
{ L_n} \ = \max_{z\in \Gamma} \sum_{k=0}^n {{|\omega_n (z)|}\over {|\omega_n^\prime (z_{n,k}) } (z-z_{n,k} ) |} 
$$
directly. Their values, together with the quotients over $log(n)$,  for $n = 2^k, k = 4,5,\ldots, 21$, are listed in Table 2.

\begin{table}[h]
\centering
\begin{tabular}{|r|r|r|}
\hline
         $n = 2^k$ &   $L_n$ & $L_n/log(n)$ \\
\hline
16	    &   4.838368	&   1.745072      \\
\hline
32	    &   5.291439	&   1.526787      \\
\hline
64	    &   6.618634	&   1.591445      \\
\hline
128	    &   8.423336	&   1.736043      \\
\hline
256	    &   11.747927	&   2.118584      \\
\hline
512	    &   12.597528	&   2.019377      \\
\hline
1024	&   14.007973	&   2.020923      \\
\hline
2048	&   15.57379	&   2.042566      \\
\hline
4096	&   17.093169	&   2.055019      \\
\hline
8192	&   18.318128	&   2.032882      \\
\hline
16384	&   20.026551	&   2.063729      \\
\hline
32768	&   21.788232	&   2.095585      \\
\hline
65536	&   24.609175	&   2.218971      \\
\hline
131072	&   26.356276   &   2.23671	      \\
\hline
262144	&   28.538463   &   2.28735	      \\
\hline
524288	&   30.692394	&   2.330514      \\
\hline
1048576	&   34.066354	&   2.457368      \\
\hline
2097152	&   36.405524   &   2.501051	 \\
\hline
\end{tabular}  
\caption{Lebesgue constants and the quotients over $log(n)$ for $n = 2^k$}
\end{table}

It is  believed that the Lebesgue constants grow to infinity at least as fast as $O(log(n))$ for large values of the polynomial degree $n$.  From our numerical experimentation, we are surprised to note that while the logarithm rate of growth is indeed $O(log(n))$ for $256\le n \le 32768$, they seem to tend to infinity with growth order $O((log(n))^2)$, for polynomial degrees $n >32768$.
\clearpage

\begin{figure} [h]
\centering
%\begin{minipage}{0.6\textwidth}
\includegraphics[width=6in,height=3in]{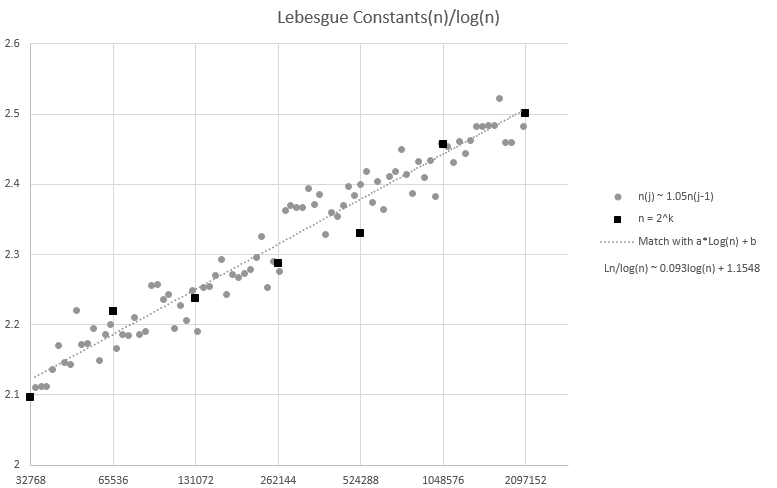} 
%\end{minipage}
\caption{Quotients of Lebesgue constants over ${log} (n)$ for even $n, n = 33994, 35694,\ldots,1950430$}
\end{figure}

Figure 3 is the graph of the quotients of $L_n$ over ${log}(n)$, for even $n$, $n = 33994, 35694,\ldots,1950430$, each number is about $105\%$ of the previous $n$. We also draw some $L_n/log(n)$ for $n = 2^k, k = 13, 14, \ldots, 20$ in the graph.  We found the small linear tendency of $L_n/{log}(n)$ under the scale of $log(n)$. With linear matching, it seems $L_n/log(n) = 0.093log(n) + 1.1548$

This numerical experimentation suggests that the Lebesgue constants could have faster growth order than smooth arcs, so that the corner in the middle of $L$-shape arc does alter the growth of the Lebesgue constants.

\subsection{Muckenhoupt $A_p$-weight constants}
Recall from Theorem 1 that the key to the validity of the Marcinkiewicz-Zygmund inequalities is the $M$-$A_p$ weight condition. In this sub-section, we report our numerical experimental findings concerning the boundedness of the Muckenhoupt $A_p$-weight constants introduced in (2.2). Let $z_0 = \psi(\rho_n e^{it_0})$ be so chosen that $|\omega(z_0)|=\min_{\zeta\in \Gamma_n}|\omega_n (\zeta)|$, where  $\Gamma_n$ is the level curve with $\rho_n = 1 + 1/(n+1)$. We then calculate 
\begin{equation}
M_n = \sup \Bigg( {1\over {|I|}} \int_I |\omega_n^p (\zeta) | d \zeta \Bigg)^{1\over p} \ \Bigg({1\over {|I|}} \int_I |\omega_n^{-q} (\zeta) | d \zeta \Bigg)^{1\over q},
\end{equation}
where the supremum is taken over all $I\subset \Gamma_n$;  and  approximate the 2 integrals in (3.5) by

$$\sum \Big| \omega^p_n (\psi (\rho_n e^{it_k}) ) \Big| \Big| \psi(\rho_n e^{it_{k+1}}) - \psi(\rho_n e^{it_{k}}) \Big|$$

and
$$\sum \Big| \omega^{-q}_n (\psi (\rho_n e^{it_k}) ) \Big| \Big| \psi(\rho_n e^{it_{k+1}}) - \psi(\rho_n e^{it_{k}}) \Big|,$$ respectively, with $t_k = t_0 + {k\over {128}}{\pi\over {n+1}}$. To validate the step size of ${\pi\over {128(n+1)}}$ is small enough, we also tried the step size of ${\pi\over {256(n+1)}}$ and obtained very similar numerical results. In Table 3, we show the Muckenhoupt $A_p$-weight constants for $p = 2,4,8,16,32,64$ and $128$, and various values of $n$.

\begin{table}[h]
\centering
\begin{tabular}{|r|r|r|r|r|r|r|r|}
\hline
   $n = 2^k$ &        $p=2$ &      $p = 4$ &      $p = 8$ &     $p = 16$ &     $p = 32$ &     $p = 64$ &    $p = 128$ \\
\hline
        16 &       1.59 &       1.66 &       1.81 &       1.95 &       2.06 &       2.14 &       2.18 \\
\hline
        32 &       2.47 &       2.41 &       2.58 &       2.76 &       2.91 &       3.02 &       3.09 \\
\hline
        64 &       2.24 &       2.30 &       2.70 &       3.07 &       3.34 &       3.51 &       3.62 \\
\hline
       128 &       2.65 &       2.69 &       3.09 &       3.49 &       3.77 &       3.95 &       4.06 \\
\hline
       256 &       3.92 &       3.57 &       3.92 &       4.31 &       4.60 &       4.79 &       4.91 \\
\hline
       512 &       4.88 &       3.93 &       4.19 &       4.59 &       4.91 &       5.12 &       5.25 \\
\hline
      1024 &       5.71 &       4.36 &       4.52 &       4.92 &       5.27 &       5.52 &       5.67 \\
\hline
      2048 &       7.73 &       5.24 &       5.24 &       5.62 &       6.01 &       6.30 &       6.47 \\
\hline
      4096 &      10.05 &       6.11 &       5.89 &       6.24 &       6.65 &       6.96 &       7.17 \\
\hline
\end{tabular}

\caption{Muckenhoupt $A_p$-weight constants on $\Gamma_n$, $n = 2^k$}
\end{table}
In Figure 4, the Muckenhoupt $A_p$-weight constants $M_n$ for $p = 2$ and are shown for even $n = 10, 12, \ldots, 4096$. In view of the divergence trend of the sequence ${M_n}$, we suspect that (2.2) is violated, in that the constants $0<c_1 \le c_2$ in (2.1) might not exist, so that the Fej\'er points (adjusted or not) would not constitute an $M$-$Z$ family for the typical $L$-Shape arc $\Gamma$ for $1<p<\infty$. This will be further studied in the next sub-section.

\begin{figure}[h]
\centering
%\begin{minipage}{0.6\textwidth}
\includegraphics[width=6in,height=3in]{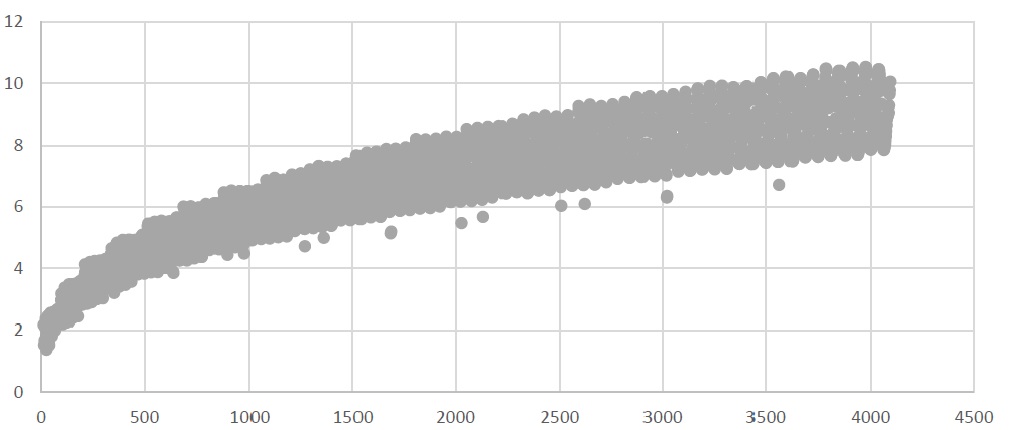} 
%\end{minipage}
\caption{Muckenhoupt $A_p$-weight constants for $p = 2$ on $\Gamma_n$ with even $n$, for $n = 10,12,\ldots,4096$}
\end{figure}
\clearpage
\subsection{Marcinkiewicz-Zygmund family}
As alluded by the observation in the previous subsection, we present our numerical experimentation in exploring if  the canonical Lagrange interpolating polynomials $P_{n,k}(z)$, defined by $P_{n,k}(z_{n,j})= \delta_{j,k}$, the Kronecker $\delta$ symbol, would satisfy (1.5), for the Fej\'er points $z_{n,k}$ to qualify as a  Marcinkiewicz-Zygmund ($M$-$Z$) family. Let
\begin{equation}
R^p_{n,k} := {{\int_{\Gamma} |P_{n,k}(z)|^p |dz|}\over{dist(z_{n,k}, \Gamma_n)}}.
\end{equation}
Then a necessary condition for the Fej\'er points $z_{n,k}$ to be an $M$-$Z$ family is the existence of some positive constant $c_2$, such that
$$R^p_{n,k} \le c_2,$$
for $k=0,\ldots,n$ and all positive integers $n$. For each $n$, let $k = k(n)$ be so chosen that $z_{n, k(n)}$ is sufficiently close to $\zeta_n\in \Gamma_n$ where $|\omega(\zeta_n)| = min_{\zeta \in\Gamma_n}|\omega(\zeta)|$, and consider 

\begin{equation}
R^p_{n} :=  R^p_{n,k(n)} = {{\int_{\Gamma} |P_{n,k(n)}(z)|^p |dz|}\over{dist(z_{n,k(n)}, \Gamma_n)}},
\end{equation}

In addition, we may also consider:
\begin{equation}
C_{p,n} := max_{0\le k \le n} {{\int_{\Gamma} |P_{n,k}(z)|^p |dz|}\over{\sum_{0\le j \le n} |P_{n,k}(z_{n,j})|^p dist(z_{n,j}, \Gamma_n)}} = max_{0\le k \le n} R^p_{n,k},
\end{equation}

In figure 5 and 6, we plot the values of $R^p_n$ for $p = 2$ and $1.5$, respectively, and for even $n = 10, 12, \ldots, 4096$.

\begin{figure} [h]
\centering
%\begin{minipage}{0.6\textwidth}
\includegraphics[width=6in,height=3in]{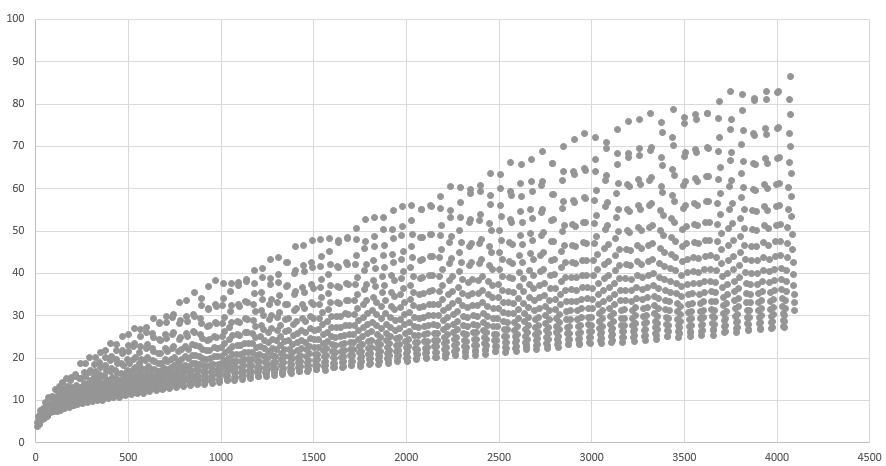} 
%\end{minipage}
\caption{$R^2_n$ for even $n, n = 10,12,\ldots,4096$}
\end{figure}

\begin{figure} [h]
\centering
%\begin{minipage}{0.6\textwidth}
\includegraphics[width=6in,height=3in]{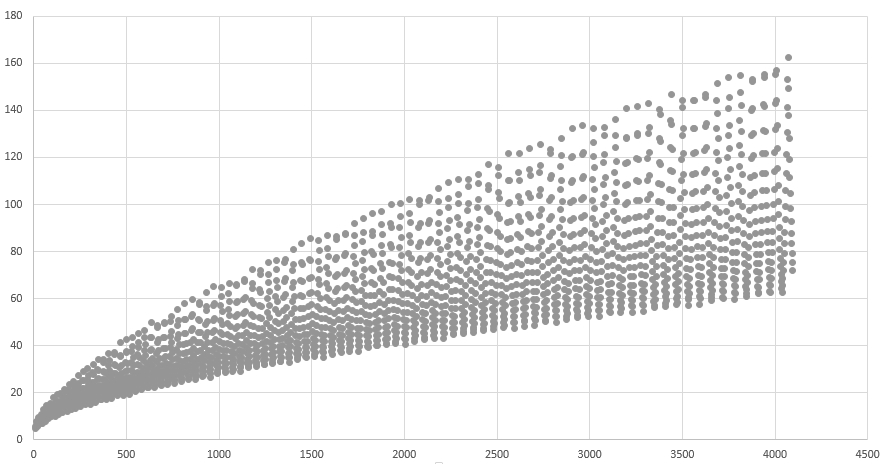} 
%\end{minipage}
\caption{$R^{1.5}_n$ for even $n, n = 10,12,\ldots,4096$}
\end{figure}

 \begin{figure} [h]
\centering
%\begin{minipage}{0.6\textwidth}
\includegraphics[width=6in,height=3in]{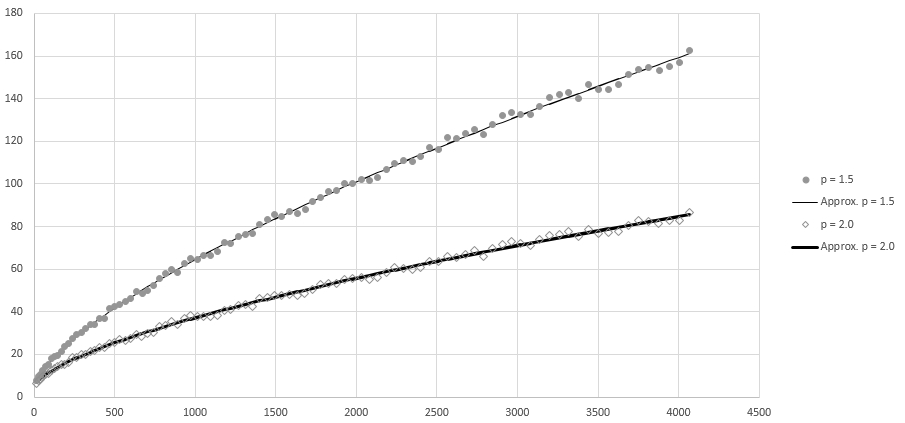} 
%\end{minipage}
\caption{Least-square curve fittings of the sub-sequences ${R_n^{1.5}}$ and ${R_n^{2}}$}
\end{figure}
\clearpage

In Figure 7, to convince ourselves that the values of ${R_n^{1.5}}$ and ${R_n^{2}}$ are unbounded, we plot the least-square approximation by the curves $a + bn^\beta$, with parameters $a, b$ and $\beta$. The two least-square approximating curves, for $p=1.5$ and $p = 2$, shown in Figure 7, are

\begin{itemize}
  \item[(a)]
$p = 1.5$: $a = 3.999147, b = 0.553288, \beta = 0.679688$,
	
 \item[(b)]
$p = 2$:    $a = 4.574969, b = 0.369888, \beta = 0.648438$,
\end{itemize}
These two monotonic increasing curves provide the trends of the values of ${R_n^{1.5}}$ and ${R_n^{2}}$, for increasing degrees of the canonical Lagrange polynomials.

\subsection{Concluding remarks}

While the Marcinkiewicz-Zygmund inequalities provide a useful tool for the assessment of suitable sampling points $\{z_{n,k}: k = 0,\cdots, n\}$ for polynomial interpolation or approximation, based on comparison of the representing polynomial $P_n$ at the data samples $P_n(z_{n,k})$  with the continuous measurement of $P_n(z)$ for $z\in \Gamma$, there have been few available criteria for selecting such sampling locations. In this paper, suitability of the sampling points on an $L$-shape arc $\Gamma$ is characterized by the $M$-$A_p$ weights based on $\omega_n(z) := \Pi_{k = 0}^n(z - z_{n,k})$ (together with certain separation condition on $\{ z_{n,k}\}$). In this regard, we believe the main result (Theorem 1), established in Section 2, should be valid for any piece-wise smooth open arc (with finitely many corners). Unfortunately, the proof of this extension seems to be much harder and is therefore delayed to a future work.

On the other hand, from numerical experimentation, it is clear that the behavior of $|\omega_n (z)|$ for the typical $L$-shape arc is significantly different from that for the  interval $[-1, 1]$ and smooth open arcs $\gamma$. This indicates that the Marcinkiewicz-Zygmund inequalities might fail to hold for the canonical Lagrange interpolating polynomials at the  Fej\'er points on $L$-shape arcs and more generally for piece-wise smooth open arc with corners, even with adjustment. The growth order of the Lebesgue constants for the adjusted Fej\'er points on $L$-shape arcs seems to be faster than that of the Chebyshev points on $[-1,1]$. Theoretic development in this direction is also delayed to our future research.

\end {document}